\documentclass{amsart}
\usepackage{amscd,amssymb,graphicx,epsfig}
\usepackage{array}

\title[Compression of Finite Groups Actions]{Compression of Finite Group Actions and 
Covariant Dimension}

\author{Hanspeter Kraft and Gerald W. Schwarz}

\address{Hanspeter Kraft\newline
\indent Mathematisches Institut der
Universit\"at Basel,\newline
    \indent Rheinsprung 21, CH-4051 Basel, Switzerland}
\email{Hanspeter.Kraft@unibas.ch}

\address{Gerald W. Schwarz \newline
    \indent Department of Mathematics\newline
    \indent Brandeis University\newline
    \indent PO Box 549110\newline
    \indent Waltham, MA 02454-9110}
\email{schwarz@brandeis.edu}
\date{September 2006}

\thanks{The first author is partially supported by the Swiss National Science
    Foundation (Schweizerischer National\-fonds),
   and the second author by  NSA Grant H98230-04-1-0070}
   
\dedicatory{To Ernest Vinberg on the occasion of his 70th birthday}
 
\newtheorem{thm}{Theorem}[section]
\newtheorem*{thm*}{Theorem}

\newtheorem{prop}{Proposition}[section]
\newtheorem{lem}{Lemma}[section]

\newtheorem{cor}{Corollary}[section]
\newtheorem*{cor*}{Corollary}
\newtheorem*{conj*}{Conjecture}

\theoremstyle{definition}
\newtheorem{defn}{Definition}[section]

\theoremstyle{remark}
\newtheorem*{rem*}{Remark}
\newtheorem{rem}{Remark}[section]

\newcommand{\op}{\operatorname}

\newcommand{\pp}{{\op{\mathfrak p}}}

\newcommand{\name}[1]{\textsc{#1\/}}
\newcommand{\NN}{{\mathbb N}} 
\newcommand{\QQ}{{\mathbb Q}}
\newcommand{\ZZ}{{\mathbb Z}}
\newcommand{\PP}{{\mathbb P}}
\newcommand{\FF}{{\mathbb F}}
\newcommand{\CC}{{\mathbb C}}

\newcommand{\Cst}{{{\mathbb C}^*}}

\newcommand{\OOO}{\mathcal O}
\newcommand{\CCC}{\mathcal C}

\newcommand{\simto}{\xrightarrow{\sim}}
\newcommand{\be}{\begin{enumerate}}
\newcommand{\ee}{\end{enumerate}}

\newcommand{\Id}{\op{Id}}

\newcommand{\SLtwo}{{\op{SL}_2}}
\newcommand{\SL}{\op{SL}}

\newcommand{\GL}{\op{GL}}

\newcommand{\PGL}{\op{PGL}}
\newcommand{\Hom}{\op{Hom}}
\newcommand{\Ind}{\op{Ind}}

\newcommand{\Aut}{\op{Aut}}

\newcommand{\Ker}{\op{Ker}}

\newcommand{\cdim}{\op{covdim}}
\newcommand{\edim}{\op{edim}}
\newcommand{\ch}{\op{char}}
\newcommand{\ps}{\par\smallskip}
\renewcommand{\Im}{\op{Im}}
\newcommand{\phimax}{{\phi_{\max}}}
\newcommand{\gr}{\op{gr}}
\newcommand{\inv}{^{-1}}
\newcommand{\pr}{\op{pr}}
\newcommand{\tphi}{{\tilde\phi}}
\newcommand{\id}{\op{Id}}
\newcommand{\quot}{/\!\!/}
\newcommand{\Span}{\op{span}}
\newcommand{\HH}{\op{H}}
\newcommand{\MM}{\op{M}}
\newcommand{\GCD}{\op{GCD}}
\newcommand{\SU}{\op{SU}}

\newcommand{\tdeg}{\op{tdeg}}
\newcommand{\Knt}{{\widetilde{K_n}}}

\renewcommand{\phi}{\varphi}

\begin{document}

\begin{abstract} 
Let $G$ be a finite group and $\phi\colon V\to W$ an equivariant morphism of finite dimensional $G$-modules. We say that $\phi$ is faithful if
$G$ acts faithfully on $\phi(V)$. The covariant dimension of $G$ is the minimum of the dimension of $\overline{\phi(V)}$ taken over all
faithful $\phi$. 

In this paper we investigate covariant dimension and are able to determine it for abelian groups and to obtain estimates for the symmetric and alternating groups. 
We also classify groups of covariant dimension less %
than 
or equal to 2. 

A byproduct of our investigations is the existence of a purely transcendental field of definition of  
degree $n-3$ for a generic field extension of degree $n\geq 5$.
\end{abstract}

\maketitle

\tableofcontents

\section{Introduction} \label{sec:introduction}

Our base field is the field $\CC$ of complex number. It could be replaced by any algebraically closed field of characteristic zero. Let $G$ be a finite group. All $G$-modules that we
consider will be finite-dimensional over $\CC $.

\begin{defn} 
A {\it covariant\/} of $G$ is an equivariant morphism $\phi\colon V \to W$ where $V$ and $W$ are
 $G$-modules. The {\it dimension\/} of $\phi$ is defined to be the dimension of the image of $\phi$:
$$
\dim \phi := \dim\overline{\phi(V)}.
$$
The covariant $\phi$ is   {\it faithful\/} if the group $G$ acts faithfully on the image $\phi(V)$. Equivalently, there is  a point $w\in\phi(V)$ with trivial isotropy group $G_w$.
\end{defn}
Here is a slightly different point of view  \cite{Rei04}.

\begin{defn} Let $V$ be a $G$-module and $X$ a faithful affine $G$-variety. A $G$-equivariant dominant morphism $\phi\colon V \to X$  is called a {\it compression}.
\end{defn}
Clearly, a faithful covariant $\phi\colon V \to W$ defines a compression $\phi\colon V \to X:=\overline{\phi(V)}$, and every compression arises in this way.
We are interested in finding compressions (faithful covariants) with small dimension. This leads to the following definition.
\begin{defn} 
The {\it covariant dimension\/} of $G$ is defined to be the minimum of $\dim X$ where $\phi\colon V \to X$ runs  over all compressions of $G$. Equivalently,
$$
\cdim G := \min\{\dim \phi \mid \phi\colon V \to W \text{ is a faithful covariant of $G$}\}. 
$$
\end{defn}

Suppose that $\phi\colon V\to W$ is a {\it rational map\/} which is $G$-equivariant. We call $\phi$ a {\it rational covariant}. Then one can define  the notion of $\phi$ being faithful and the dimension of 
$\phi$ as in the case of ordinary covariants.

\begin{defn} (\name{Buhler-Reichstein} \cite{BuR97}) The {\it essential dimension $\edim G$\/} of $G$ is the minimum dimension of all the faithful rational covariants   of $G$.
\end{defn}

The   covariant dimension and essential dimension of $G$  differ by at most $1$ (\cite{Rei04} and Proposition \ref{reiprop}). Essential dimension and covariant dimension of $G$ are related to cohomological invariants, generic polynomials and other topics, see
\cite{BuR97}.

In this paper we study the notion of covariant dimension staying in the category of morphisms. A major role is played by covariants $\phi\colon V\to W$ which are {\it homogeneous\/}. We are able to determine the covariant dimension of abelian groups and obtain estimates of the covariant dimension of the symmetric and alternating groups. Finally, we are able to classify the groups of covariant dimension less  %
than or equal to 2. It turns out that these are exactly the finite subgroups of $\GL_2(\CC)$. We also obtain a new result about the ``fields of definition'' for generic extensions of degree $n$. Except for this and the classification above, most of our results could also be obtained from \cite{BuR97} and \cite{Rei04}. However, we think that our methods are of independent interest.

{\small\noindent
{\bf Acknowledgement:}
The authors thank Zinovy Reichstein for introducing us to to the notions of essential dimension and covariant dimension.   We  thank Daniel Goldstein for informing us about \cite{Ga54} and his help with the cohomology groups in the last section.}

\vskip1cm
\section{First Properties} \label{sec:firstproperties}
A covariant $\phi\colon V \to W$ is called {\it minimal\/} if $\phi$ is faithful and $\dim\phi = \cdim G$. 
\begin{lem}\label{lem-covariant}
Let $V,W$ be  two $G$-modules and let $v\in V, w\in W$ be such that $G_v \subset G_w$. Then there is a covariant $\phi\colon V \to W$ such that $\phi(v) = w$.
\end{lem}

\begin{proof} By assumption  there is a $G$-equivariant  map
$\mu\colon Gv \to W$ which sends $v$ to $w$. This map lifts to a morphism from $V$ to $W$  which we can average
over $G$ to obtain   a covariant $\tilde \mu$ extending $\mu$.
\end{proof}

\begin{rem}\label{Rem1}  Obviously, one can prescribe the images $w_1,w_2,\dots,w_m\in W$ of a finite number of points
$v_1,v_2,\dots,v_m\in V$ from distinct orbits  provided that $G_{v_i} \subset G_{w_i}$ for all $i$.
\end{rem}

\begin{prop} Let $V,W$ be two faithful  $G$-modules and let $v\in V$ and $w\in W$ be points with trivial stabilizer. Then
there is a minimal covariant $\phi\colon V\to W$ such that $\phi(v) = w$.
\end{prop}

\begin{proof} Let $\phi_0\colon V_0 \to W_0$ be a faithful covariant with $\dim \phi_0 = \cdim G$. Then there is a $v_0\in  V_0$ such that  $w_0:=\phi_0(v_0)\in W_0$ has a trivial stabilizer. Thus $v_0$ has a trivial stabilizer, too. By the previous lemma we can find   covariants $\phi_1\colon V \to V_0$ and $\phi_2\colon W_0 \to W$ and points $v\in V$ and $w\in W$ with trivial stabilizer such that 
$\phi_1(v) = v_0$ and
$\phi_2(w_0) = w$.  Then
$\phi:= \phi_2\circ\phi_0\circ\phi_1$ is faithful with $\dim\phi\leq\cdim G$, hence we have equality. 
\end{proof}

Here are some elementary properties of covariant dimension. We leave the proofs to the reader.

\begin{rem}\label{remgeneral} 
(a) Let $H$ be a subgroup of $G$. Then $\cdim H\leq \cdim G$. 

(b) If $G$ is a product $G_1\times G_2$, then $\cdim
G\leq\cdim G_1+\cdim G_2$.

(c) If $G$ is non-trivial cyclic, then $\cdim G=1$.
\end{rem}
Moreover, Remark~\ref{remgeneral} holds for  essential dimension in place of covariant dimension. From \cite{Rei04} we have

\begin{prop} \label{reiprop} Let $G$ be a finite group. Then $\edim G\leq \cdim G\leq\edim G+1$.
\end{prop}

\begin{proof} The first inequality is clear. Let $\phi\colon V\to W$ be a rational faithful covariant of minimal  dimension. Then there is 
a nonzero $f\in\CC [V]^G$ such that $\Phi:=f\phi$ is a covariant. 
Now $\overline{\Phi(V)}$ is contained in the cone on $\overline{\Im\phi}$, so $\cdim G\leq\edim G+1$.
\end{proof}

\vskip1cm
\section{Covariant Dimension for Abelian Groups} \label{sec:abeliangroups}

Let $G$ be a finite abelian group. We can write $G = G_1\times G_2\times\cdots \times G_n$ where $G_i$ is cyclic of order $d_i$ and $d_1|d_2|\cdots | d_n$. Then $n$ is the {\it rank\/} of $G$.  In this section we show that $\cdim G=n$. From Remark~\ref{remgeneral} we have 
$\cdim G\leq n$.

\begin{lem}\label{lemabelian}
Let $p$ be a prime number and let $K$ be a field of characteristic $0$ or $p$. Let $f_i\in
K[x_1,x_2,\ldots,x_n]$, $i=1,\ldots,n$ be polynomials of the form
\begin{align*}
f_1(x_1,x_2,\ldots,x_n) &\in K[x_1,x_2^p,x_3^p,\ldots,x_n^p]\setminus K[x_1^p,x_2^p,\ldots,x_n^p]\\
f_2(x_1,x_2,\ldots,x_n) &\in K[x_1^p,x_2,x_3^p,\ldots,x_n^p]\setminus K[x_1^p,x_2^p,\ldots,x_n^p]\\
&\vdots\\
f_n(x_1,x_2,\ldots,x_n) &\in K[x_1^p,x_2^p,x_3^p,\ldots,x_n]\setminus 
K[x_1^p,x_2^p,\ldots,x_n^p]
\end{align*}
Then the  Jacobian determinant   $\det (\frac{\partial f_i}{\partial x_j})$ is non-zero. In particular, $f_1,f_2,\ldots,f_n$ are algebraically independent.
\end{lem}

\begin{proof} (a) If $\ch K = p>0$,  then 
$(\frac{\partial f_i}{\partial x_j})$ is a diagonal matrix with non-zero entries $\frac{\partial f_i}{\partial x_i}$, and the 
lemma follows.
\ps
(b) If $\ch K = 0$ we use a ``reduction mod $p$'' argument to reduce to case (a). Let $C\subset K$ be the set of coefficients of  the polynomials $f_i$ and set $L:= \QQ(C)$. We  can find algebraically independent elements $a_1,\ldots a_m\in L$ such that the elements of  $C$ are algebraic over $\QQ(a_1,\ldots,a_m)$. By multiplying the polynomials $f_i$ with suitable elements from $L$ we can assume that the elements of $C$ are integral over $\ZZ[a_1,\ldots,a_m]$. Thus $R := \ZZ[a_1,\ldots,a_m][C]$ is a free $\ZZ$-module and $f_i\in R[x_1,\ldots,x_n]$. Then $pR \subsetneqq R$, and we can assume that $f_i \not\equiv 0 \mod p$. Now it
follows from (a) that $\det (\frac{\partial f_i}{\partial x_j})\not\equiv 0 \mod p$, 
hence the lemma.
\end{proof}

\begin{rem} Lemma~\ref{lemabelian} does not hold if the characteristic of $K$ is positive and prime to
$p$. In fact,   $\det(\frac{\partial f_i}{\partial x_j})$ vanishes for $p=2$, $f_1 := x_1^3, f_2:=x_2^3\in\FF_3[x_1,x_2]$. Of course, $f_1,f_2$ are still algebraically independent, but we do not know if this holds in general.
\end{rem}

Recall our decomposition $G=G_1\times G_2\times\dots\times G_n$ where $d:=d_1|d_2|\dots|d_n$ and $d_i=|G_i|$. Fix  embeddings $G_i \subset \CC ^*$. The homomorphism 
$$
g=(\zeta_1,\zeta_2,\ldots,\zeta_n) \mapsto \bmatrix \zeta_1   & & &\\ & \zeta_2  & &\\ &&\ddots &\\ 
&&& \zeta_n 
\endbmatrix \in \GL_n(\CC)
$$
defines a faithful  representation of $G$ of dimension $n$. 
\begin{thm}\label{thm-abelian}
Let $G$ be an abelian group of rank $n$ and $V$ a faithful representation of dimension $n$. Then any faithful covariant $\phi\colon V\to V$ is dominant. Hence $\cdim G=n$. 
\end{thm}

\begin{proof} 
It is enough to prove the theorem for the faithful representation defined above.
Fix a prime divisor $p$ of $d$.   It suffices to show that the components $\phi_i$  of  
$\phi$ are of the form given in Lemma 2. Write $\phi_1 =
\sum_{k=0}^m h_k(x_2,\ldots,x_n) x_1^k$. Then the $h_k$'s are invariants for the subgroup $G_2\times \cdots\times G_n$, hence
$h_k\in \CC[x_2^p,\ldots,x_n^p]$. On the other hand, $\phi_1$ is a covariant for $G_1$ and so
$\phi_1(\zeta_1x_1,x_2,\ldots,x_n) = \zeta_1 \phi_1(x_1,x_2,\ldots,x_n)$ %
for $\zeta_1\in \CC$, which implies that $h_k = 0$ unless
$k \equiv 1\mod p$. Hence we have that $\phi_1$ is of the form $f_1(x_1^p,x_2^p,\ldots,x_n^p) x_1 $, and similarly for $\phi_2,\dots, \phi_n$.
\end{proof}

\begin{rem} The theorem holds for $\edim G$ in place of $\cdim G$ \cite{BuR97}.
\end{rem}

\vskip1cm
\section{Faithful Groups and Irreducible Covariants} \label{sec-irreduciblecov}

We investigate conditions under which we may assume that
$\cdim G$ is realized by  an ``irreducible'' and homogeneous covariant. We start with the following easy lemma. 

\begin{lem}\label{lem-faithful} Let $W=\oplus_{i=1}^r W_i$ be  faithful where the $W_i$ are irreducible. 
Let $\phi=(\phi_1,\dots,\phi_r)\colon V \to W$ be a covariant. If $\phi_i\neq 0$, $i=1,\dots,r$, then 
$\phi$ is faithful.
\end{lem}

\begin{proof} 
If $\phi$ is not faithful, then $N:=\Ker(G\to\Aut(\overline{ \phi(V))})$ is a nontrivial normal subgroup of $G$. For some $i$ we must have that $W_i^N=\{0\}$, so that $\phi_i=0$.
\end{proof}

\begin{defn} We say that the covariant $\phi\colon V\to W$ is {\it irreducible\/} if $W$ is irreducible. We say that the group $G$ is {\it faithful\/} if it has a faithful  irreducible module.
\end{defn}

Note that the symmetric group $S_n$ is faithful as is any product  of
simple groups. Also, if $G$ is faithful, then by Schur's Lemma, the center of $G$ must be a cyclic group. In general, there is the following useful criterion for a group $G$ to be faithful, due to \name{Gasch\"utz}.

\begin{prop}[\cite{Ga54}]\label{prop-Gaschutz}
 Let $G$ be a finite group and denote by $N_G \subset G$ the subgroup generated by the  minimal normal abelian subgroups. Then $G$ is faithful if and only if $N_G$ is generated by the conjugacy class of one of its elements.
\end{prop} 

\begin{cor}\label{cor-Gaschutz}
Let $G$ be a non-faithful group and $H \subset G$ a subgroup containing $N_G$. Then $H$ is non-faithful, too.
\end{cor}

\begin{proof} Since $N_G \subset N_H$ are both products of cyclic groups of prime order we can write $N_H = N_G \times M$ with        a suitable normal subgroup $M\subset H$. By Proposition~\ref{prop-Gaschutz}, $N_G$ can't be generated by the $G$-conjugacy class of a single element and so  $N_G\times M$ can't be generated by the $H$-conjugacy class of a single element.
\end{proof}

Let $\phi=\sum_{j\leq n} \phi_j\colon V\to W$ be a covariant where $\phi_j$ is
homogeneous of degree $j$, $1\leq j\leq n$. Assume that $\phi_n$ is not identically
zero. We call $\phi_n$ the {\it maximal homogeneous component\/} and denote it by $\phimax$.

\begin{lem}\label{lem-phimax} Let $\phi\colon V \to W$ be as above. Then $\dim\phimax\leq \dim\phi$.
\end{lem}

\begin{proof} Let $X$ denote $\overline{\phi(V)}$, and let $\pp$ denote the (prime) ideal
of $X$ in $\OOO(W)$. If $f\in\pp$, let $\gr f$ denote the highest degree nonzero
homogeneous part of $f$, and let $\gr \pp$ be the ideal generated by all the $\gr f$,
$f\in \pp$. Then $\CCC  X$, the associated cone of $X$, is the zero set of $\gr \pp$, and
$\dim X=\dim \CCC  X$ (see \cite{Kr85}). We show that $\Im \phimax\subset \CCC  X$, which
gives the lemma.

Let $f=\sum_{i=0}^m f_m$ be in $\pp$, where $\gr f=f_m\neq 0$. Then for $v\in V$,
$0\neq t\in \CC$ we have $0=f(\phi(t\inv v))$ which implies that
$$
f_m(\phi(t\inv v))=-\sum_{j=0}^{m-1} f_j(\phi(t\inv v)).
$$
Multiplying both sides by $t^{nm}$, %
where $n=\deg \phi$, we obtain
$$
f_m(t^n\phi(t\inv v))=-\sum_{j=0}^{m-1} t^{n(m-j)} f_j(t^n\phi(t\inv v)).
$$
Letting $t$ go to zero we see that the left hand side above converges to
$f_m(\phimax(v))$ and the right hand side converges to zero. Thus $f_m$ vanishes on
$\Im \phimax$, i.e., $\Im\phimax\subset \CCC X$.
\end{proof}

If $\phi$ is faithful, it is not clear that $\phimax$ is also. However, for faithful groups this is almost automatic. In fact, we have the following result which is an immediate consequence of Lemmas~\ref{lem-faithful} and \ref{lem-phimax} above.

\begin{prop}\label{prop-faithful} 
Let $G$ be faithful with irreducible faithful $G$-module $W$ and faithful $G$-module $V$. Then there is a homogeneous minimal covariant $\phi\colon V\to W$. 
\end{prop}

\begin{rem} Starting with a faithful representation $V$ we cannot always guarantee that there is a faithful homogeneous covariant $\phi\colon V \to V$ with minimal dimension. Let $V=\CC^2$ where  $\ZZ/2$   (resp.\ $\ZZ/3$) acts by multiplication by roots of unity on the first (resp.\ second) copy of $\CC$. Let $x$ and $y$ be coordinate functions on $V$. Then a faithful minimal covariant for $G:=\ZZ/2\times\ZZ/3$ is $\phi(x,y)=(x^3y^3,x^4y^4)$. Suppose that we had a homogeneous
faithful covariant $\psi$. Then for some $(x_0,y_0)\in V$, $\psi(x_0,y_0)=(x_1,y_1)$ where $x_1y_1\neq 0$. By equivariance, $\psi(-x_0,y_0)=(-x_1,y_1)$. Thus  the image of $\psi$, which is a cone, contains two linearly independent vectors. It follows that  $\overline{\Im
\psi}=V$ and so $\dim\psi=2> \cdim G = 1$.
\end{rem}

\begin{cor}\label{cor-edim}
Let $G$ be a faithful group with trivial center. Then 
$$
\edim G = \cdim G -1.
$$
\end{cor}
\begin{proof} Let $\phi\colon V \to V$ be a homogeneous minimal covariant. Then the rational covariant $\psi\colon V \to V \to \PP(V)$ is faithful of dimension $\cdim G - 1$.
\end{proof}

The usefulness of the existence of {\it homogeneous\/} minimal covariants for calculating the covariant dimension is shown by  
Proposition~\ref{prop-addfactor} below. We first need a definition.
For a set $X$, let $|X|$ denote its cardinality and let $\Id_X$ denote the identity map of $X$.

\begin{defn} Let $V$ be a $G$-module and
$\rho_V\colon G\to\GL(V)$ the corresponding representation. Define 
$$
z_G(V)=z(V):= |\rho_V(G) \cap\Cst\Id_V|.
$$

If $V$ is irreducible and faithful, we have $z_G(V) = |Z(G)|$.
\end{defn}

\begin{prop}\label{prop-addfactor}
Assume that $G$ has a homogeneous minimal covariant $\phi\colon V \to V$. If $m>0$ is coprime to $z_G(V)$, then $\cdim G\times \ZZ/m = \cdim G$.
\end{prop}

For the proof we will use the following result.

\begin{lem}\label{lem-deginvariants} 
There is an integer $n_0>0$ and an open dense set $V'\subset V$ with the following property: For every $n\geq n_0$ there is a homogeneous invariant $h\in\OOO(V)^G$ of degree $n\cdot z_G(V)$ which has no zeroes on $V'$.
\end{lem}

\begin{proof} Set $d:=z_G(V)$.
Let $M \subset \NN$ be the monoid of degrees  
of homogeneous    
elements of $\OOO(V)^G$. 
By assumption, we have $M \subset d\NN$, and so the subgroup $\langle M \rangle\subset  \ZZ$ generated by $M$ equals  $d'\ZZ$  
for some multiple $d'$ of $d$. It follows that all $G$-invariants are also invariant under $\mu_{d'}\subset \Cst$, the $d'$th roots of unity. Thus $\mu_{d'} \subset \rho_V(G)$ and so $d'=d$.

Since $\langle M \rangle = d\ZZ$ we can find two homogeneous invariants $h_1,h_2\in \OOO(V)^G$ with $\gcd(\deg h_1,\deg h_2) =  
d$. Therefore, for every $n$ large enough there is a monomial $h_1^\alpha h_2^\beta$ of degree $n\cdot  
d$. The  
lemma now follows by setting $V'$ to be the complement of the zero set of $h_1h_2$.
\end{proof}

\begin{proof}[Proof of Proposition~\ref{prop-addfactor}] Let $f\in\OOO(V)^G$ be a non-zero homogeneous invariant. Then 
$f\cdot \phi\colon V \to V$ is again a homogeneous faithful covariant of minimal dimension. In fact,  there is a $v\in V$ such that $f(v)\neq 0$ and such that  $\phi(v)$ has a trivial stabilizer. Then the same holds for any non-zero multiple 
$\lambda\phi(v)$ and, in particular, for $f(v)\phi(v)$. Thus $f\cdot\phi$ is faithful. Since $\phi$ is homogeneous the image $\phi(V)$ is a cone and so $(f\cdot\phi)(V) \subset \phi(V)$, hence $f\cdot\phi$ is of minimal dimension.

The covariant $f\cdot\phi$ has degree $\deg\phi+\deg f$ which, by Lemma~\ref{lem-deginvariants},  can be any number of the form $\deg\phi + n\cdot z_G(V)$ for $n\geq n_0$.  Since $m$ is coprime to $z_G(V)$ there is an $n\geq n_0$ such that 
$\deg\phi + n\cdot z_G(V)$ is divisible by $m$. For the corresponding covariant $f\cdot\phi$  this implies that it is also equivariant with respect to the scalar action of  $\ZZ/m$ on $V$, and so $f\cdot\phi$ is a covariant for $G\times\ZZ/m$.
\end{proof}

\begin{cor}\label{cor-addZm}
If $G$ is faithful and $m>0$ coprime to $|Z(G)|$, then $\cdim G \times \ZZ/m = \cdim G$. %
Moreover, $G\times \ZZ/m$ is faithful.
\end{cor}

\begin{cor}\label{cor-semidirect}
Let $G = \ZZ/3^\ell \ltimes A$, $\ell\geq 1$,  where $A$ is a finite abelian $2$-group of rank 2 and a generator $a$ of $\ZZ/3^\ell$ acts non-trivially on $A$. 
Then $\cdim G \geq 3$.
\end{cor}

\begin{proof}  
Denote by $\alpha$ the automorphism of $A$ induced by $a$. We can assume that $\alpha$ is trivial on $2A$
since otherwise 
$2A$ has again rank 2 and we can
replace $G$ by the subgroup $\ZZ/3^\ell\ltimes 2A$.  
Then the induced automorphism on $A/2A$ is non-trivial since the order of $\alpha$ is not a power of $2$. It follows that the automorphism $\alpha$ is given by a matrix of the form
$\bmatrix 2r & 1+2s \\ 1+2t & 1+2u \endbmatrix$ with respect to suitable generators of $A$. Since $\alpha$ is the identity on $2A$ we get $2A=0$ and so $A \simeq \ZZ/2\oplus\ZZ/2$. 
It follows from Gasch\"utz's Criterion (Proposition~\ref{prop-Gaschutz}) that $G$ is faithful: $N_G = \langle a^3,(\ZZ/2)^2\rangle$ and is generated by the conjugates of $a^3\cdot(1,0)$. The center $Z(G)$ is generated by $a^3$ and so $\cdim G = \cdim G\times \ZZ/2 \geq \cdim (\ZZ/2)^3 = 3$.
\end{proof}

\vskip1cm
\section{Existence of Homogeneous Covariants} \label{sec:homogeneouscov}

If $v\in V$ is a {\it principal\/} point, i.e., the stabilizer of $v$ equals the kernel of $\rho_V\colon G \to \GL(V)$, and $w\in V$ is arbitrary, we can always find a covariant $\phi\colon V \to V$ such that $\phi(v) = w$  (see Lemma~\ref{lem-covariant}). In order to find a {\it homogeneous\/} covariant with this property, we need an additional assumption.

\begin{prop}\label{prop-homcovariant}
Let $V$ be a $G$-module and let $v\in V$ be a principal point such that the corresponding point $[v]\in \PP(V)$ is also principal. Let $w$ be another point of $V$. Then there is a homogeneous covariant 
$\phi\colon V\to V$ such that $\phi(v)=w$. 
\end{prop}

For the proof we need the following result which is probably well-known.

\begin{lem}\label{lem-vandermonde}
Let $V$ be a vector space of dimension $\geq 2$ and $v_1,v_2,\ldots,v_s \in V$ pairwise linearly independent elements. If $r\geq s-1$, then $v_1^r,\ldots,v_s^r \in S^rV$ are linearly independent.
\end{lem}
 
\begin{proof} (a) We first consider the case $\dim V = 2$. By choosing a suitable basis and multiplying the $v_i$'s with scalars if necessary we can assume that $v_i = (1,b_i)$, $i=1,\ldots ,s$. Then $v_i^r = (1,b_i, b_i^2, \ldots, b_i^r)\in S^r(V)\simeq \CC^{r+1}$ and  
linear independence follows from the nonsingularity of  the Vandermonde matrix $(b_i^j)$.

(b) There is alway a linear projection $\rho\colon V\to W$, $\dim W = 2$, such that the images $\rho(v_1),\ldots,\rho(v_s)$ are pairwise linearly independent.  
So the general case follows from part (a).
\end{proof}
          
\begin{proof}[Proof of Proposition~\ref{prop-homcovariant}]
Define $H:=\{g\in G \mid gv \text{ is a scalar multiple of } v\}$. By assumption, $H = \Ker (G \to \PGL(V))$.
We have a character $\chi$ of $H$, where $h\cdot v'=\chi(h)v'$ for any $v'\in V$. Let
$g_1 H,\dots, g_s H$ be the set of left cosets of $H$ where $g_1=e$. Then,  by Lemma~\ref{lem-vandermonde}, the elements $g_i\cdot v$ give rise to linearly independent elements $g_i \cdot v^r$ of $S^rV$, for $r \geq s-1$. Choose $r$ to be congruent to $1$ modulo $d:=|H|$. Then the submodule $W:=\Span\{v^r,g_2\cdot v^r,\ldots,g_s\cdot v^r\}\subset S^rV$ is isomorphic to  the representation $\Ind_H^G \CC_\chi$ induced  from the character $\chi$ of $H$. Now there is a linear $H$-equivariant map $\CC_\chi \to V$ such that $1\mapsto w$,  and so the induced map 
$W\simto \Ind_H^G\CC_\chi \to V$ is $G$-equivariant and sends $v^r$ to $w$. It follows that the composition 
\begin{equation*}
\begin{CD}
V @>{v\mapsto v^r}>> S^rV @>{\pr}>> 
W \simto
\Ind_H^G \CC_\chi @>>> V
\end{CD}
\end{equation*}
is the required homogeneous covariant, where $\pr$ is equivariant projection onto $W$.
\end{proof}

\begin{rem} One can establish a more general form of the proposition. For $v\in V$ let $\tilde G_v:=\{g\in G\mid g\cdot v=\lambda v$ for some $\lambda\}$. There is the obvious character $\chi_v\colon\tilde G_v\to\Cst$. Let $W$ be another $G$-module and $w\in W$. Then the following are equivalent:
\begin{enumerate}
\item $\tilde G_v\subset \tilde G_w$ and $\chi_w(g)=\chi_v(g)^d$ for all $g\in\tilde G_v$ and for some $d\in\NN$.
\item There is a homogeneous covariant $\phi\colon V\to W$ such that $\phi(v)=w$.
\end{enumerate}
\end{rem}
Now we can prove the existence of a homogeneous minimal covariant under some additional assumption on the faithful representation $V$ without assuming that $V$ is irreducible. We first need a lemma about the degree of a covariant.

\begin{lem}\label{lem-degcovariant}
If $\phi\colon V \to V$ is a covariant, then $\deg \phi \equiv 1 \mod z_G(V)$. 
\end{lem}

\begin{proof} Since every homogeneous component of $\phi$ is a covariant we can assume that $\phi$ is homogeneous, say of degree $d$. Then $\phi(t\cdot v) = t^d\cdot \phi(v)$ for all $t\in \CC$, $v\in V$. Now choose $g\in G$ so that $\rho_V(g) = \zeta \id_V$ where 
$\zeta$ is a primitive $z(V)$th root of unity. Then we have 
$$
\phi(gv) = \phi (\zeta \cdot v) = \zeta^d \cdot \phi(v)\quad \text{and}\quad g\phi(v) = \zeta \cdot \phi(v)
$$
for all $v\in V$. Hence, $\zeta^d = \zeta$ which implies that $d\equiv 1 \mod z(V)$, as claimed.
\end{proof}

\begin{prop}\label{prop-homminimal}
Let  $V$  be a faithful representation and let $V = \bigoplus_{i=1}^n V_i$  be a decomposition into irreducible submodules. Assume that $z(V_i) = z(V)$  for $i<n$ and that every prime divisor of $z(V_n)$ divides $z(V)$. Then there is a homogeneous minimal covariant $\phi\colon V \to V$. 
\end{prop}
\begin{proof} We give the proof for $n=2$ and leave the obvious generalization to the reader.
Let $\phi\colon V \to V$ be a faithful minimal covariant. We can clearly assume that the two components $\phi_1,\phi_2$ are both non-zero. If $\deg \phi_1 = \deg\phi_2$ then we are done: $\phimax$ has two non-zero components ${\phi_1}_{\max}$ and ${\phi_2}_{\max}$, hence is faithful (Lemma \ref{lem-faithful}), and $\dim\phimax \leq \dim\phi$ by Lemma~\ref{lem-phimax}.
\par\smallskip
We reduce to the case above by composing $\phi$ with a covariant $\psi$ of the form 
\begin{equation}\label{equ1}
\psi(v_1,v_2) = (f_1(v_1)v_1,f_2(v_2)v_2)= (\psi_1(v_1), \psi_2(v_2))
\end{equation}
where $f_1\in\OOO(V_1)^G$ and $f_2\in\OOO(V_2)^G$ are homogeneous invariant functions, so that the two components of  the composition  $\phi\circ\psi$ are both non-zero and have the same degree. 
 
It follows from Lemma~\ref{lem-deginvariants} applied to the two representations $V_1$ and $V_2$ that there are open dense subsets $V_1'\subset V_1$, $V_2'\subset V_2$ and an integer $n_0>0$ such that, for every $n\geq n_0$ there are homogeneous invariants $f_i\in\OOO(V_i)^G$ of degree $nz(V_i)$ which have no zeroes in $V_i'$, ($i=1,2$). 

Since $V_i$ is irreducible, the image of ${\phi_i}_{\max}\colon V \to V_i$ contains a principal point $v_i$ such that $[v_i]\in\PP(V_i)$ is also principal. By Proposition~\ref{prop-homcovariant}, there is a homogeneous covariant $\mu_i\colon V_i \to V_i$ such that $\mu_i(v_i) \in V_i'$. Replacing $\phi_i$ by the composition $\mu_i\circ\phi_i$ we can therefore assume that the image of ${\phi_i}_{\max}$ meets $V_i'$ $(i=1,2)$. Then 
the compositions $\psi_i\circ{\phi_i}_{\max}$ are non-zero if the invariants $f_i$ in equation~(\ref{equ1}) are chosen according to Lemma~\ref{lem-deginvariants}.   
Set $\tphi_i:=\psi_i\circ\phi_i$, $i=1$, $2$.  
We have
$$
\deg \tphi_1 = \deg \phi_1 \cdot (1+\deg f_1) \quad
\text{and}
\quad 
\deg \tphi_2 = \deg \phi_2 \cdot (1+\deg f_2).
$$
Now $\deg\phi_i = 1 + a_i \cdot z(V)$ (Lemma~\ref{lem-degcovariant}), and so we have to solve the equation
$$
(1+x z(V_1))(1+ a_1z(V)) = (1+y z(V_2))(1+ a_2z(V))
$$ 
with integers $x$, $y \geq n_0$ which is possible by the following lemma. 
\end{proof}

\begin{lem}\label{lem-solvequation} 
Let $d,d_1,d_2,a_1,a_2 \in \NN$. Assume that $d|d_1|d_2$ and that $d_2$ has the same prime divisors as $d$ (i.e., $d_2|d^N$ for large $N$). Then the equation 
\begin{equation*}
(1+a_1d)(1+xd_1) = (1+a_2d)(1+yd_2)
\end{equation*}
has a solution $x$, $y\in \NN$ if and only if $a_1\equiv a_2 \mod\frac{d_1}{d}$. Moreover, $x$ and $y$ can be chosen to be arbitrarily large.
\end{lem}
\begin{proof}
The equation implies   
that $a_1d \equiv a_2d \mod d_1$, hence  $a_1\equiv a_2 \mod\frac{d_1}{d}$. Conversely,  assume that $a_1d \equiv a_2d \mod d_1$. Then $1+a_1d \equiv 1+ a_2d \mod d_1$. Moreover, $1+a_1d$ is invertible mod $d_2$ and so there is an $m\in \NN$ such that $(1+ a_1d)m \equiv 1 \mod d_2$. It follows that $(1+a_2d)m\equiv 1 \mod d_1$, hence
$$
(1+a_1d)m = 1+yd_2 \quad\text{and}\quad (1+a_2d)m = 1+xd_1
$$
for some $x$, $y\in\NN$, and so $(1+a_1d)(1+xd_1) = (1+a_2d)(1+yd_2)$. The last statement is clear.
\end{proof}

\begin{cor}\label{cor-Z2M}
Let $M$ be an abelian group of rank $k$ whose order is odd, and let $G:=\ZZ/2^s\ltimes M$, $s\geq 1$, be a  semidirect product where a generator $a$ of $\ZZ/2^s$ acts on $M$ by sending each element to its inverse.  
Then $\cdim G \geq k+1$.
\end{cor}
\begin{proof} Replacing $G$ by a suitable subgroup we can assume that $M=(\ZZ/p)^k$.
The group $\bar G := \ZZ/2^s \ltimes \ZZ/p$ has a two-dimensional faithful irreducible representation. Using the $k$ copies of $\ZZ/p$ we therefore obtain $k$ irreducible representations $V_1, \cdots V_k$ of $G$ such that $V := \bigoplus_iV_i$ is faithful. The center $\bar Z$ of $\bar G$ is generated by $a^2$ and is isomorphic to $\ZZ/2^{s-1}$. Thus, the assumptions of Proposition~\ref{prop-homminimal} are satisfied with $z(V)=z(V_i)=2^{s-1}$, and we can find a homogeneous minimal covariant $\phi\colon V \to V$. Hence, by Corollary~\ref{cor-addZm} and Theorem~\ref{thm-abelian}, $\cdim G = \cdim G\times \ZZ/p \geq \cdim (\ZZ/p)^{k+1} = k+1$.
\end{proof}

\vskip1cm
\section{Products with Cyclic Groups} \label{sec:productwithcyclic}

We have seen in section~\ref{sec-irreduciblecov} that $\cdim G \times \ZZ/p = \cdim G$ if $G$ is faithful and $p$ is coprime to $|Z(G)|$. We will show now that $\cdim G\times \ZZ/p = \cdim G + 1$ if $p$ is a divisor of $|Z(G)|$. More generally, we have the following.
\begin{prop}\label{prop-addH}
Let $V=W\oplus \CC_\chi$ be a faithful representation of $G$ where $W$ is  irreducible and $\chi$ is a character of $G$. Assume that $z(W)$ and $|\chi(G)|$ have the same prime divisors as  $z(V)$ and that either $z(W)=z(V)$ or $|\chi(G)|=z(V)$. If the kernel $H$ of the action of $G$ on $W$ is non-trivial, then 
$\cdim G = \cdim G/H + 1$.
\end{prop}

\begin{proof} By assumption, we have an embedding $G \hookrightarrow G/H \times \chi(G)$. Since $\chi(G)$ is cyclic we get $\cdim G \leq \cdim (G/H\times\chi(G))\leq \cdim G/H + 1$.

For the reverse inequality, we consider a faithful covariant $\phi\colon V \to V$ of minimal dimension. By Proposition~\ref{prop-homminimal} we can 
assume that $\phi$ is homogeneous. Set $m:=|\chi(G)|$. We have $\phi(w,t) = (F(w,t^m), t\cdot h(w,t^m))$ where $h$ is a $G$-invariant function. If $F(w,0)=h(w,0) = 0$, then, because $\phi$ is homogeneous, we can divide it by $t^m$ without changing its dimension 
or  faithfulness. Thus we can assume that either $F(w,0)\neq 0$ or $h(w,0)\neq 0$. In the first case, $\phi|_W \colon W \to W$ is non-zero,  
hence a faithful covariant of $G/H$. Since $\phi(W) \subset \phi(V)^H \subsetneq \phi(V)$ we see that $\dim \overline{\phi(V)} > \dim \overline{\phi(W)}\geq \cdim G/H$,  
as desired.

If $h(w,0)\neq 0$ and $F(w,t^m) = t^{sm} F_0(w,t^m)$, $s\geq 1$, where $F_0(w,0)\neq 0$ we define 
$$
\psi(w,t):=
(\frac{F(w,t)}{(t\cdot h(w,t^m))^{sm}},t\cdot h(w,t^m)) = (\frac{F_0(w,t)}{h(w,t^m)^{sm}},t\cdot h(w,t^m)).
$$ 
The equivariant morphism $\psi=(\psi_1,\psi_2)$ is defined on the dense open set $V_{h}\subset V$ where $h$ does not vanish. Moreover, $\dim \overline{\psi(V_h)} = \dim \overline{\phi(V)}=\cdim G$, because  $\psi_1$ and $\psi_2$  generate the same subfield of $\CC(V)$ as the two components of $\phi$. By definition, $\psi_1(w,0)\neq 0$, and $\psi(W_h) = \psi_1(W_h) = F_0(W_h)$  
since $\psi_1$ is homogeneous of non-zero degree. Since $F_0|_W\colon W \to W$ is a non-zero (hence faithful) covariant for $G/H$, we get $\dim\overline{\psi(W_h)} \geq \cdim G/H$. Finally, $\overline{\psi(W_h)} \subset \overline{\psi(V_h)}^H \subsetneq \overline{\psi(V_h)}$, and so
$$
\cdim G = \dim\overline{\psi(V_h)}> \dim\overline{\psi(W_h)} \geq \cdim G/H.
$$
\end{proof}

\begin{rem}\label{rem-proof}
The proof above has two parts. First one shows that there is a homogeneous minimal covariant $\phi\colon V \to V$ for $V = W \oplus \CC_\chi$, and then one proves the inequality $\cdim G \geq \cdim G/H + 1$. The assumptions about $z_G(V), z_G(W)$ and $|\chi(G)|$ are only used in the first part. Once the existence of a homogeneous minimal covariant $\phi$  is established, then the proof above applies
if the action of $G$ on $W$ has a non-trivial kernel $H$.
\end{rem}

\begin{cor}\label{cor-center}
Let $G$ be a faithful group and $p$ a prime divisor of  $|Z(G)|$. Then 
$\cdim G\times\ZZ/p = \cdim G + 1$.
\end{cor}
\begin{proof} This follows from the proposition above by choosing for $W$ an irreducible faithful representation of $G$ (with trivial action of $\ZZ/p$) and for $\chi$ the standard character $r\mapsto e^{2\pi i \frac{r}{p}}$ of $\ZZ/p$.
\end{proof}

\begin{cor}\label{cor-twoZ/p}
If $H$ is a faithful group and $q$ a prime which does not divide $|Z(H)|$, then $\cdim H\times (\ZZ/q)^2 = \cdim H + 1$.
\end{cor}
\begin{proof}
This is clear from Corollary~\ref{cor-center}: Take $G = H\times \ZZ/q$. Then $G$ is faithful and $\cdim G = \cdim H$.
\end{proof}

We finish this section  
with a result on the covariant dimension of the semi-direct product of two cyclic $2$-groups. In the proof we will need a modification of the proof of Proposition~\ref{prop-homminimal} in order to reduce to a homogeneous minimal covariant.

\begin{prop} \label{prop-Z2Z2}
Let $G = \ZZ/2^k \ltimes \ZZ/2^\ell$, $k,\ell \geq 1$. If $G$ is commutative or $k=1$, then $\cdim G = 2$. 
Otherwise, $\cdim G \geq 3$.
\end{prop}
%
%
\begin{proof} (a) We can assume that $G$ is not commutative and so $\ell\geq 2$. 
Let $a$ be a generator of $\ZZ/2^k$, let $b$ be a generator of $\ZZ/2^\ell$, and denote by $\alpha$ the (non-trivial) automorphism induced by $a$ on $\ZZ/2^\ell$. If $k=1$, then $\alpha(b) = b^{-1}$ if $l=2$ and $\alpha(b)=b^{-1}$, $b^{2^{l-1}+ 1}$ or $b^{2^{l-1}- 1}$ if $l>2$, and one easily constructs a faithful representation on $\CC^2$, see part (b). Thus we may assume that $k\geq 2$. We now show that $\cdim G\geq 3$.

(b) We may assume that $\alpha$ acts trivially on $2\ZZ/2^\ell$, and then it follows that $\alpha$ has order $2$ and sends $b$ to $b^{2^{l-1}+1}$.  
Let $\xi$ be a primitive $2^\ell$th root of unity and let $\tau$ be a primitive $2^k$th root of unity, We have the irreducible representation $\rho_W$ of $G$ on $W:=\CC^2$ where $b$ acts by the diagonal matrix $\bmatrix \xi&0 \\ 0&\xi^{2^{l-1}+1}\endbmatrix$ and   $a$ acts  
by the matrix $\bmatrix 0 & 1\\  1 & 0\endbmatrix$. 
We also have the character $\chi\colon G \to \Cst$ defined by $\chi(a) := \tau$  and $\chi(b) :=1$. The direct sum $V:=W\oplus\CC_\chi$ is a faithful representation $\rho$ of $G$, so, clearly, $\cdim G\leq 3$. We have $z_G(W) = 2^{\ell -1}$ since $\rho_W(b^2)$ is a scalar matrix. 

If $k\geq \ell$  we can assume that  $\xi=\tau^{2^{k-\ell}}$. Then $\rho(b^2 a^{2^{k-\ell+1}}) = \xi^2\Id_V$ and so $z_G(V) = 2^{\ell-1} = z_G(W)$. Thus we can apply Proposition~\ref{prop-addH} and find that $\cdim G = \cdim \rho_W(G) + 1 = 3$.
\par\smallskip
There remains the case where $k<l$. Then we set $\tau=\xi^{2^{\ell-k}}$ and find that the scalar matrices 
in $\rho(G)$ are generated by $\rho(a^2b^{2^{\ell-k+1}})$. Thus $z_G(V) = 2^{k-1}<z_G(\CC_{\chi})=2^k\leq  z_G(W)=2^{\ell -1}$ and we cannot apply Proposition~\ref{prop-addH} directly, but need a modification.

Let $\phi=(F,h)\colon V\to V$ be a minimal covariant where $F\colon V\to W$ and $h\colon V\to\CC_\chi$ are $G$-equivariant. If $\deg F = \deg h$, then $\phi_{\max} = (F_{\max}, h_{\max})$ is again faithful and minimal, and  
we may proceed
as in Proposition~\ref{prop-addH} (see Remark~\ref{rem-proof}). 

In general,  $\deg F$ and $\deg h$ are both  $\equiv 1 \mod 2^{k-1}$, by Lemma~\ref{lem-degcovariant}. Consider the two homogeneous invariants $f_1(x,y):=xy^{2^{\ell-1}-1} + yx^{2^{\ell-1}-1} \in \OOO(W)^G$ and $f_2(t):= t^{2^k}\in\OOO(\CC_\chi)^G$ of degree $2^{\ell-1}$ and $2^k$, respectively, and the corresponding covariants
$$
\psi = ({f_1}^u\Id_W, {f_2}^v\Id_{\CC_\chi})\colon V \to V
$$
where $u$, $v\in\NN$. 
Clearly, $h_{\max}$ and $h$ are dominant, as are $F_{\max}$ and $F$ (else we get that $\cdim \rho_W(G)\leq 1$). 
Thus the composition $\tilde\phi:=\psi\circ\phi$ is faithful and  the two components $\tilde F = {f_1}(F)^u\cdot F$ and $\tilde h = f_2(h)^v \cdot h$ have degrees
$$
\deg \tilde F = (1+ u2^{\ell-1})\deg F\quad
\text{and}\quad
\deg \tilde h = (1+v2^{k})\deg h.
$$
Setting  $\deg F = (1+ r2^{k-1})$ and $\deg h = (1+ s2^{k-1})$ it suffices to solve the equation
\begin{equation}\label{equ2}
(1+ r2^{k-1})(1+ u2^{\ell-1}) = (1 + s 2^{k-1})(1+ v2^{\ell-1}).
\end{equation}
By Lemma~\ref{lem-solvequation} this is possible if and only if $r$ and $s$ have the same parity, i.e., if and only if $\deg F \equiv \deg h \mod 2^k$. (In the notation of that lemma we have 
$d=2^{k-1}, d_1=2^k, d_2=2^{\ell-1}$.)

If $\deg F \not\equiv \deg h \mod 2^k$, consider the following 
$G$-invariant function on $V$. 
$$
 f(x,y,t):= (xy^{2^{\ell-1}-1}-x^{2^{\ell-1}-1}y)t^{2^{k-1}}.
$$
Then $f$ is bihomogeneous of degree $(2^{\ell-1}, 2^{k-1})$.
Since $h_{\max}\neq 0$ and $F_{\max}$ is dominant, 
$f(F_{\max},h_{\max})$ is non-zero and 
$\tilde F := f(F,h)\cdot F\colon V \to W$ has degree 
$$
\deg\tilde F = (2^{\ell-1}+1)\deg F + 2^{k-1} \deg h \equiv (1+(r+1)2^{k-1}) \mod 2^k .
$$ 
As a consequence,  $\tilde\phi := (\tilde F , h)$ is a minimal covariant and $\deg\tilde F \equiv \deg h\mod 2^k$. Now we can apply Lemma~\ref{lem-solvequation} and solve equation~(\ref{equ2}). This finishes the proof.
\end{proof}

\vskip1cm
\section{Covariant Dimension and Essential Dimension for $S_n$ and $A_n$} \label{sec:Sn}

If we put together our results so far for the faithful groups $S_n$ $(n\geq 3)$  and $A_n$ $(n\geq 4)$ we have the following.

\begin{thm}\label{thm-Sn} 
Let $n\geq 3$.
\begin{enumerate}
\item $\cdim S_n=\edim S_n+1$.
\item $\cdim S_n\times\ZZ/d=\cdim S_n$ for $d\geq 2$.
\item $\cdim S_n\geq \lfloor\frac n2\rfloor+1$.
\item $\cdim S_{n+2}\geq\cdim S_n+1$.
\ee
\end{thm}

\begin{thm}\label{thm-An} 
Let $n\geq 4$.
\begin{enumerate}
\item $\cdim A_n=\edim A_n+1$.
\item $\cdim A_n\times\ZZ/d=\cdim A_n$, for $d\geq 2$.
\item $\cdim A_n\geq 2 \lfloor\frac n4\rfloor+1$.
\item $\cdim A_{n+2}\geq\cdim A_n+1$.
\ee
\end{thm}

\begin{proof}[Proof of Theorem~\ref{thm-Sn} and \ref{thm-An}] 
Parts (a) and (b) follow from Corollary~\ref{cor-edim} and \ref{cor-addZm}.
For part (c) we use (b) and Theorem~\ref{thm-abelian}. Since
$S_{2m}\supset (\ZZ/2)^m$, we have $\cdim S_n=
\cdim S_n\times\ZZ/2\geq\cdim ((\ZZ/2)^{\lfloor\frac n2\rfloor+1})=\lfloor\frac n2\rfloor+1$. One proceeds similarly for $A_n$ using that $A_{4m}\supset (\ZZ/2)^{2m}$. 

Finally, for part (d)  we have $\cdim S_{n+2}=\cdim S_{n+2}
\times\ZZ/2\geq\cdim (S_n\times (\ZZ/2)^2)\geq \cdim
(S_n\times\ZZ/2)+1=  (\cdim S_n)+1$, by Corollary~\ref{cor-center}, and similarly for $A_n$
\end{proof}

These results allow us to determine the covariant dimension for the small symmetric and alternating groups. Just recall that $\cdim A_n \leq \cdim S_n \leq n-1$ since there  is a
faithful $n-1$-dimensional representation. 
$$
\cdim S_ 2 = \cdim A_3 = 1,\quad \cdim S_3 =2, \quad \cdim A_4 = \cdim S_4 = 3.
$$
We will see below that $\cdim S_n \leq n-2$ for $n\geq 5$, hence 
$$
\cdim A_5 = \cdim S_5 = 3 \quad\text{and}\quad  \cdim S_6 = 4.
$$ 
The first unknown cases are $\cdim A_6$ which is either $3$ or $4$, and  $\cdim S_7$ 
which is either $4$ or $5$.

\medskip
\begin{tabular}{ | c || c | c | c | c | c | c |}
\hline
$n$ & \makebox[1cm]{2} &  \makebox[1cm]{3} &  \makebox[1cm]{4} & \makebox[1cm]{5} &  \makebox[1cm]{6} & \makebox[1cm]{7}\\
\hline
$\cdim S_n$  & 1 & 2 & 3 & 3 & 4 & 4 or 5 \\
\hline
$\cdim A_n$ &  & 1 & 3 & 3 & 3 or 4 & \\
\hline
\end{tabular}

\vskip1cm
\section{Upper Bounds for the Covariant Dimension} \label{sec:upperbounds}

Let $G$ be finite group and $V$ a $G$-module such that the $G$-action normalizes a reductive subgroup $H$ of $\GL(V)$. Then we can form the {\it algebraic quotient\/} 
$$
\pi\colon V \to V\quot H
$$ 
which has coordinate ring $\OOO(V\quot H):=\OOO(V)^H \subset \OOO(V)$ (see \cite{Kr85}).
Since $G$ normalizes $H$ the invariant ring is $G$-stable and the quotient morphism $\pi$ is $G$-equivariant, hence a compression. 

\begin{lem}\label{lem-quot}
If the action of $G$ on $V\quot H$ is faithful, then
$$
\cdim G\leq\dim V\quot H \leq \dim V-\max\{\dim Hv\mid v\in V\}.
$$
\end{lem}

\begin{proof} We can find a finite dimensional $G$-stable subspace $W$ of $\OOO(V)^H$ which generates $\OOO(V)^H$. The associated
morphism $\phi\colon V\to W^*$ has image isomorphic to the quotient
$V\quot  H$. Since $\phi$ maps $V$ onto $V\quot H$ and since $G$ acts faithfully on $V\quot H$, the covariant $\phi$ is faithful. The fibers of $\phi$ have dimension  $\geq\max\{\dim Hv\mid v\in V\}$, so that $\cdim G\leq\dim V\quot H\leq \dim V-\max\{\dim Hv\mid v\in V\}$.
\end{proof}

\begin{prop}\label{prop:upperbound} 
For $n\geq 5$ we have $\cdim A_n \leq \cdim S_n \leq n-2$.
\end{prop}

For the proof we use the following construction which will also play a central r\^ole in 
section~\ref{sec:fieldext} below.

Start with the standard representation of $\SLtwo$ on $\CC^2$, and let $T'\subset \SL_n$ be the group of diagonal matrices 
of determinant $1$. Then $H:= \SLtwo\times T'$ acts linearly on $V_n:=\CC^2\otimes \CC^n$. There is also the standard action of $S_n$ by permutations on $\CC^n$, hence on $V_n$ which normalizes the action of $T'$ and commutes with the action of $\SLtwo$. 

We may 
regard an element of $V_n$ as an $n$-tuple
of elements in $\CC^2$, so there is  a canonical surjective morphism $\mu\colon V_n
\to S^n(\CC^2)$ given by multiplying the elements of $\CC^2$. This morphism is the quotient by the group $S_n\cdot T'$, and is equivariant with respect to $\SLtwo$.

Consider the quotient $\pi\colon V_n \to X_n:= V_n\quot H$. Since $S_n$ normalizes $H$ it acts on the quotient $X_n$ and $\pi$ is $S_n$-equivariant. By construction, we have a canonical isomorphism $X_n/S_n \simto S^n(\CC^2)\quot\SLtwo$. Moreover, the quotient morphism $\pi$ can be decomposed into  $\eta\colon V_n \to Y_n:= V_n\quot \SLtwo$, the quotient by $\SLtwo$, and  $\rho\colon Y_n \to Y_n\quot T'=X_n$, the quotient by $T'$, so that we obtain the following commutative diagram:
$$
\begin{CD}
V_n @>{\quot S_n\cdot T'}>{\mu}> S^n(\CC^2) \\
@V\eta V{\quot \SL_2}V @V
\pi''
V{\quot \SL_2}V  \\
Y_n @>{\quot S_n\cdot T'}>> S^n(\CC^2)\quot \SL_2\\
@V\rho V{\quot T'}V @| \\
X_n @>{ /S_n}>{\pi'}>  S^n(\CC^2)\quot \SL_2 \\
\end{CD}
$$
The generic $T'$-orbit on $V_n$ is closed with trivial stabilizer, and for $n\geq 3$ the generic $\SL_2$-orbit on $S^n(\CC^2)$ is closed with finite stabilizer. Thus, for $n\geq 3$, the generic $H$-orbit on $V_n$  is closed and one easily sees that it has a trivial stabilizer. It follows that
 $$
\dim X_n =  \dim V - \max\{\dim Hv\mid v\in V\} = 2n - (n-1+3) = n-2.
$$
The following proposition collects 
properties of the 
morphisms and quotient maps above. 

\begin{prop} \label{prop:diagram}
Consider the representation $V_n := \CC^2\otimes\CC^n$ of $H:= \SLtwo\times T'$ where $T' \subset \SL_n$ 
is the subgroup of diagonal matrices, and let $\pi\colon V_n \to X_n:=V_n\quot H$ be the quotient.
\begin{enumerate}
\item The natural action of $S_n$ on $V_n$ normalizes $H$, the quotient $\pi$ is $S_n$-equivariant and $\dim X_n = n-2$ for $n\geq 3$.

\item The quotient of $V_n$ by $\SL_2$ is given by the map $\eta\colon f_1\otimes v_1 + f_2\otimes v_2 \mapsto v_1\wedge v_2 \in \bigwedge^2\CC^n$ with image $Y_n = \GL_n (e_1\wedge e_2) \cup \{0\}$.

\item \name{(Kempe)} The invariant ring $I_n:=\OOO(Y_n)^{T'}$ is generated by the invariants of degree $d$ where $d=n$ if $n$ is odd and $d=\frac{n}{2}$ if $n$ is even. Thus the quotient map $\rho\colon Y_n \to X_n$ is homogeneous of degree $d$.

\item  The action of $S_n$ on $X_n$ and on $\PP(X_n)$  is faithful for $n\geq 5$.

\end{enumerate}
\end{prop}

\begin{proof}[Proof of Proposition~\ref{prop:upperbound}] The claim follows from Lemma~\ref{lem-quot} applied to the $S_n$-equivariant quotient $\pi\colon V_n \to X_n= V_n\quot H$, using Proposition~\ref{prop:diagram}(a) and (d).
\end{proof}

\begin{rem}\label{rem-quotSn} 
It is easy to analyze the cases $n=2$, $3$ and $4$. For $n=2$ the quotient $V_2\quot H$ is $\CC$ with a non-trivial action of  $S_2$. For $n=3$ the quotient $V_3\quot H$ is also $\CC$, but $A_2$ acts trivially. Finally, for $n=4$ the quotient $V_4\quot H$ is a hypersurface in $\CC^3$ and the representation of $S_4$ on $\CC^3$ has kernel the Klein $4$-group and quotient group $S_3$ which acts on $\CC^3$ in the standard way. 
\end{rem}

\begin{proof}[Proof of Proposition~\ref{prop:diagram}] 
Part (a) was already proved above.
Part (b) is the classical First Fundamental Theorem for $\SL_2$, see  \cite[5.2.1 Proposition and Remarks]{How88}, and
 Part (c) is due to \name{Kempe}. A proof can be found in \cite[5.4.2.5 Theorem, page 156ff.]{How88}.

For (d) it suffices to show that $S_n$ acts faithfully on $X_n$ since the only non-trivial normal subgroup of $S_n$ is $A_n$ for $n\geq 5$. But otherwise, the subgroup $A_n\subset S_n$ would act trivially on $X_n$. This means that a generic $H$-orbit on $V_n$ is stable under $A_n$ which implies that there is a non-trivial homomorphism $A_n\to H$,  a contradiction.
\end{proof}

The idea behind the proof of the upper bounds above can be used to explicitly calculate 
compressions for the group $S_n$. We have seen in Proposition~\ref{prop:diagram}(b) that the algebraic quotient $\eta\colon 
\CC^2 \otimes \CC^n \to (\CC^2\otimes \CC^n)\quot \SLtwo$ is given by 
$$
\eta\colon f_1\otimes v_1 + f_2 \otimes v_2 \mapsto v_1\wedge v_2 \in \textstyle\bigwedge^2 \CC^n
$$
and that $(\CC^2\otimes \CC^n)\quot \SLtwo \subset \bigwedge^2 \CC^n$ 
is $\{v_1\wedge v_2\mid v_1,v_2\in \CC^n\}$, the closure of the highest weight orbit in $\bigwedge^2(\CC^n)$.
Let $x_{ij}=(e_i\wedge e_j)^*$ be the usual dual basis of  
$\bigwedge^2\CC^n$. 
A monomial $\prod_{i<j} x_{ij}^{\alpha_{ij}}$ is $T'$-invariant if and only if each index $i$ occurs the same number of times, i.e., for every fixed $k$ the sum $\sum_{i<k} \alpha_{ik} + \sum_{i>k} \alpha_{ki}$ is independent of $k$. Choosing such a monomial $f$ we obtain an $S_n$-equivariant  and $T'$-invariant  morphism $\Sigma_f\colon \bigwedge^2 \CC^n \to \CC[S_n]$ given by $\sum_{\sigma\in S_n} (\sigma f)\cdot \sigma$.  We can do slightly better by composing with
the $S_n$-equivariant embedding $\iota\colon \CC^n \subset \CC^2\otimes \CC^n$, $a \mapsto f_1\otimes a + f_2\otimes (1,1,\ldots,1)$. Then $\eta\circ\iota(a_1,\ldots,a_n) = \sum_{i<j}(a_i-a_j)e_i\wedge e_j$. Now the next proposition follows immediately from what we have said so far.

\begin{prop}\label{prop-construction}
Let $f \in \CC[x_1,\ldots,x_n]$ be a monomial in the differences $(x_i-x_j)$, where $i<j$, such that each $x_i$ occurs the same number of times. If $n\geq 5$ and $f$ is not an $S_n$-invariant, then the morphism $\Sigma_f\colon \CC^n \to \CC[S_n]$  
corresponding to   $\sum_{\sigma\in S_n} (\sigma f)\cdot\sigma$ defines a compression of dimension $\leq n-2$.
\end{prop}

\begin{rem} If $n$ is even, we can use $f := (x_1-x_2)(x_3-x_4)\cdots(x_{n-1}-x_n)$ and get a compression of dimension $\leq n-2$ and of degree $\frac{n}{2}$. In general, we can always use $f := (x_1-x_2)(x_2-x_3)\cdots (x_{n-1}-x_n)(x_n-x_1)$ to obtain a compression of dimension $\leq n-2$ and of degree $n$. As a consequence of \name{Kempe}'s result (Proposition~\ref{prop:diagram}(c)) one shows  that these compressions have dimension equal to $n-2$.

It is an open problem if there exist other monomials $f$ such that the corresponding covariant $\Sigma_f$ has dimension strictly less than $n-2$. So far, all our explicit calculations have only produced covariants of dimension $n-2$.
\end{rem}

\vskip1cm
\section{Purely transcendental  fields of definition for generic field extensions} \label{sec:fieldext}
In this section we shortly describe the relation between the essential dimension of a generic field extension of degree $n$ and the essential dimension of $S_n$ due to \name{Buhler-Reichstein}, and then show that every such field extension is defined over a purely transcendental extension of degree $n-3$.

In the following we assume that all fields contain the complex numbers $\CC$. Let $L/K$ %
be a finite field extension.

\begin{defn} We say that $L/K$ is {\it defined over a subfield $K' \subset K$\/} if there is %
a finite field extension  $L'/K'$ of degree $[L':K'] = [L:K]$ such that $L = L'\, K$. The minimal transcendence degree (over $\CC$) of such a subfield $K'$ is called the {\it essential dimension\/} of $L/K$:
$$
\edim (L/K) := \min \{\tdeg_\CC K' \mid \text{$L/K$ is defined over $K'$}\}.
$$
\end{defn}
Now  assume that $\tdeg_\CC K = \infty$ and  consider  the general field extention 
$K_n/K$ of degree $n$ defined by the equation
\begin{equation}\label{eq1}
x^n + a_1 x^{n-1} + a_2 x^{n-2} + \cdots + a_{n-1} x + a_n = 0
\end{equation}
where the $a_i \in K$ are algebraically independent. The following result is due to \name{Buhler-Reichstein}, \cite[Corollary~4.2]{BuR97}. It was the starting point for studying compressions of group actions.

\begin{thm}[\name{Buhler-Reichstein}] $\edim K_n/K = \edim S_n$.
\end{thm}
In order to prove the inequality  $\edim K_n/K \leq \edim S_n$  one shows that every rational faithful $S_n$-covariant $\rho\colon \CC^n \to W$ of dimension $d$ determines a subfield $K' \subset K$ of $\tdeg_\CC K' = d$ such that $K_n/K$ is defined 
over $K'$. We will need this construction in a slightly different form which we are going to explain now.

Denote by $\Knt \supset K_n$ the splitting field of the equation (\ref{eq1}) so that $\Knt/K$ is a Galois extension with Galois group $S_n$. Clearly, $\Knt$ contains elements $x_1,x_2,\ldots,x_n$ which are permuted under $S_n$ and generate $\Knt/K$. By assumption, the elements $x_1,\ldots,x_n$ are algebraically independent over $\CC$. Define $V := \CC x_1 + \cdots + \CC x_n \subset \Knt$. By construction, $V$ is the standard representation of $S_n$.  

Now let $\phi\colon V \to W$ be a homogeneous $S_n$-covariant of dimension $d$ and consider the cone $X := \overline{\phi(V)} \subset W$. The field $\CC(X)$ of rational functions on $X$ can be considered as a subfield of $\CC(V) = \CC(x_1,\ldots,x_n)$ and $\tdeg_{\CC} \CC(X) = d$. Moreover, the field $\CC(X)$ contains  the subfield $\CC(\PP(X))$ of rational functions on the projective variety $\PP(X):= X \setminus\{0\}/\Cst$, i.e., the subfield generated by all quotients $\frac{p}{q}$ where $p,q$ are homogeneous regular functions on 
$X$ of the same degree.

Assume now that the $S_n$-action on $\PP(X)$ is faithful (which is always the case if $\phi$ is faithful and $n\geq 5$). Then $\CC(\PP(X))/\CC(\PP(X))^{S_n}$ is a Galois extension with Galois group $S_n$ and so $K \cdot \CC(\PP(X))^{S_n} = \Knt$. This shows that the extension $\Knt/K$ is defined over $\CC(\PP(X))^{S_n}$ and the same holds for $K_n/K$. Thus we have proved the following result.

\begin{prop}\label{prop:defining-field} 
Let $K_n/K$ be the field extension defined by the equation
\begin{equation*}
x^n + a_1 x^{n-1} + a_2 x^{n-2} + \cdots + a_{n-1} x + a_n = 0
\end{equation*}
where the coefficients $a_1,a_2,\ldots,a_n\in K$ are algebraically independent over $\CC$, and let $\phi\colon \CC^n \to W$ be a homogeneous covariant. Define $X := \overline{\phi(\CC^n)}$ and assume that $S_n$ acts faithfully on $\PP(X)$. Then $K_n/K$ is defined over a subfield isomorphic to  $\CC(\PP(X))^{S_n}$.
\end{prop}

In order to apply this result we use the explicit construction of a faithful covariant for $S_n$ given in 
section~\ref{sec:upperbounds}, using the representation of $H:= \SLtwo\times T'$ on $V_n:=\CC^2\otimes \CC^n$ together with the action of $S_n$ on $V_n$ by permutations normalizing $H$. 
As before, we embed $\CC^n$ into $V_n=\CC^2\otimes\CC^n$ by $a \mapsto f_1\otimes a + f_2 \otimes (1,1,\ldots,1)$, and obtain an $S_n$-equivariant linear map $\iota\colon \CC^n \to V_n$. The composition $\eta\circ\iota$ is the linear map $a=(a_1,\ldots,a_n) \mapsto \sum_{i<j}(a_i-a_j)\,e_i\wedge e_j$ whose kernel is the trivial representation $\CC \subset \CC^n$.

\begin{prop} \label{prop:explicit-covariant}
\begin{enumerate}
\item The composition $\phi:=\rho\circ\eta\circ\iota \colon \CC^n \to X_n$ is a homogeneous covariant of degree $d$ where $d=n$ if $n$ is odd and  $d=\frac{n}{2}$ if $n$ is even. 
\item $\phi\colon \CC^n \to X_n$ is surjective. 
\item For $n\geq 5$, the action of $S_n$ on $\PP(X_n)$ is faithful and $\dim \PP(X_n) = n-3$. 
\item The varieties $\PP(X_n)$ and $\PP(X_n)/S_n$ are both rational, for all $n$.
\end{enumerate}
\end{prop}
\begin{proof}  
Part (a) follows from Proposition~\ref{prop:diagram}(c). For (b) it suffices to show that $\pi'\circ\phi\colon \CC^n \to 
S^n(\CC^2) \quot\SL_2$ is surjective since $\pi'\colon X_n \to 
S^n(\CC^2) \quot\SL_2$ is the quotient by the finite group $S_n$. Now the composition $\gamma\colon \CC^n \to V_n \to S^n(\CC^2)$ is given by $a\mapsto \prod_{i}(a_ix+y)$, and so the image in $S^n(\CC^2)$ meets every $\SLtwo$-orbit except $\{0\}$. Thus, $\pi'\circ\phi=
\pi''\circ\gamma$ is surjective.


(c) This follows from Proposition~\ref{prop:diagram}(a) and (d).

(d) The field $\CC(\PP(X_n))^{S_n}$ is isomorphic to 
$\CC(S^n(\CC^2))^{\SLtwo\times\Cst}$ and $\CC(\PP(X_n))$ is isomorphic to $\CC(V_n)^{\SLtwo\times T'\times\Cst}$, and the latter fields are both rational, due to a result of  \name{Katsylo} \cite{Kat84}.
\end{proof}

Combining this result with Proposition~\ref{prop:defining-field} we  get the following consequence.

\begin{thm} Let $K_n/K$ be the field extension of degree $n\geq 5$ defined by the equation
\begin{equation*}
x^n + a_1 x^{n-1} + a_2 x^{n-2} + \cdots + a_{n-1} x + a_n = 0
\end{equation*}
where the coefficients $a_1,a_2,\ldots,a_n\in K$ are algebraically independent over $\CC$. Then $K_n/K$ is defined over a purely transcendental extension $K'$ of $\CC$ of transcendence degree $n-3$. Moreover, if  $L'/K'$ is of degree $n$ such that $L'K = K_n$ 
then $L'/\CC$ is purely transcendental, too.
\end{thm}

\begin{rem} The two cases $n=5$ and $n=6$ are due to \name{Hermite} and \name{Joubert}, respectively. They showed that one can always find a generator of the field 
extension $K_n/K$ whose equation has the form
$$
t^5 + a t^3 + b t + b = 0 \quad 
\text{or}
\quad t^6 + a t^4 + b t^2 + c t + c =0,
$$
see \cite{Kra06}. It is unknown if similar results hold in degree $n>6$.
\end{rem}

\vskip1cm
\section{Groups of Low Covariant Dimension} \label{sec:covdim2}
In this section we describe the finite groups of covariant dimension $\leq 2$.

\begin{thm}\label{thm-covdim1}
If $G$ is a group of covariant dimension 1 then $G$ is cyclic.
\end{thm}

\begin{proof}
Let $\phi\colon V\to X$ be a one-dimensional compression. We can clearly assume that $X$ is normal, hence a rational curve, by L\"uroth's Theorem.  It follows that $X$ must be isomorphic to $\CC$ and so $G$ is a subgroup of $\Aut(\CC)$ fixing the point $\phi(0)$.  
Hence $G$ is cyclic.
\end{proof}

\begin{prop}\label{prop-faithfulcovdim2}
 If $G$ is a faithful group of covariant dimension 2, then $Z(G)$ is cyclic and $G/Z(G)$ is isomorphic to $D_{2n}$ $(n\geq 2$), $A_4$, $S_4$ or $A_5$.
\end{prop}
\begin{proof} Let $\phi\colon V \to V$ be a homogeneous faithful covariant of dimension 2 where $V$ is irreducible. Then  $X:= \overline{\phi(V)}\subset V$ is a cone of dimension 2, and so $\PP(X)\subset \PP(V)$ is an irreducible rational curve. It follows that the kernel of the action of $G$ on $\PP(X)$ equals the center $Z(G)$ of $G$. The normalization of $\PP(X)$ is the projective line $\PP^1$ and so  $G/Z(G)$ is isomorphic to a finite subgroup of $\PGL_2$. Since $G$ is not commutative the factor group $G/Z(G)$ is not cyclic and the  theorem follows.
\end{proof}

Here is the main theorem of this section;
\begin{thm}\label{thm-covdim2} Let $G$ be a group of covariant dimension 2. Then $G$ is isomorphic to a subgroup of $\GL_2$.
\end{thm}
We remark that a similar result does not hold for covariant dimension $>2$ as shown by the group $S_5$ which has covariant dimension  3.

First we show:
 
\begin{thm}\label{thm-nonfaithfulcovdim2}
If $G$ is a non-faithful group of covariant dimension 2, then $G$ is abelian of rank 2.
\end{thm}
The proof needs some preparation.

\begin{lem}\label{lem-surjective}
Let $G$ be a group of covariant dimension 2 and $V$ a faithful representation. If $a,b\in G$ 
do not commute, then there is an irreducible subrepresentation $W\subset V$ such that  the commutator  $(a,b)$ is non-trivial on $W$ and such that the image $\bar G$ of $G$ in $\GL(W)$ is a faithful group of covariant dimension 2. In particular, there is a surjective homomorphism from $G$ to $D_{2n}$ $(n\geq 2$), $A_4$, $S_4$ or $A_5$.
\end{lem}

\begin{proof} Let $\phi\colon V \to V$ be a minimal covariant and $V = \bigoplus_i V_i$ the decomposition into irreducible subrepresentations. Then the commutator $(a,b)$ acts non-trivially on the image $\phi(V)$, hence non-trivially on $\phi_i(V)$ for some component $\phi_i\colon V \to V_i$ of $\phi$. If $\bar a, \bar b$ denote the images of $a,b$ in $\GL(V_i)$, then $(\bar a,\bar b)\neq 1$ and so $\bar G\subset \GL(V_i)$, %
the image of $G$ in $\GL(V_i)$,  is non-commutative.  If $N\subset G$ is the kernel of $G\to \bar G$, then $\phi_i(V^N)\neq (0)$, and so $\phi_i\colon V^N \to V_i$ is a faithful covariant of $\bar G\simeq G/N$ of dimension $\leq \dim\phi = 2$. 
Since $\bar G$ is not commutative,  $\cdim\bar G=2$. The last statement now follows from Proposition~\ref{prop-faithfulcovdim2} above.
\end{proof}

\begin{proof}[Proof of Theorem~\ref{thm-nonfaithfulcovdim2}]
Assume that the assertion does not hold. Let $G$ be a non-commutative non-faithful group of covariant dimension 2 and of minimal order with these properties. This implies that every strict subgroup $H$ of $G$ is either faithful or commutative. In particular, if $H$ contains $N_G$ then $H$ is not faithful by Gasch\"utz's Criterion (Corollary~\ref{cor-Gaschutz}), and so $H$ is commutative. Since $G$ is not faithful with $\cdim G=2$,  we have rank $N_G=2$.
 
\smallskip
\noindent
{\bf Claim 1:}
{\it There are no surjective homomorphisms from $G$ to $A_4$, $S_4$ or $A_5$.}
\par\smallskip
If $\rho$ is a surjective homomorphism from $G$ to $A_5$ then $\rho(N_G)$ is trivial. If $\rho$  is a surjective homomorphism from $G$ to $S_4$ then $\rho(N_G)\subset K$ where $K\subset S_4$ is the Klein 4-group. In both cases $\rho^{-1}(A_4)\subsetneq G$ is  
neither faithful nor commutative, contradicting the minimality assumption.

Now assume that there is a surjective homomorphism $\rho\colon G \to A_4$, and
let $g_3\in G$ be the preimage of an element of $A_4$ of order $3$. We  
may assume that the order of $g_3$ is a power $3^\ell$. 
Then $\rho(N_G) \subset K$ and so the strict subgroup  
$S:=\rho^{-1}(K)\subsetneq G$ is commutative. Denote by $S_2$ the 2-torsion of $S$. Since $\rho(S_2) = K$ we see that $S_2$ has rank 2. Moreover, $S_2$ is normalized by $g_3$, but not centralized, and so $\cdim \langle g_3,S_2 \rangle \geq 3$ by  Corollary~\ref{cor-semidirect}. This contradiction proves Claim~1.

\smallskip
\noindent
{\bf Claim 2:} 
{\it For every prime $p>2$ the $p$-Sylow-subgroup $G_p\subset G$ is normal and commutative of rank $\leq 2$. Hence $G$ is a semidirect product $G_2\ltimes G'$  where $G':=\prod_{p>2}G_p$ and $G_2$ is a $2$-Sylow subgroup.}
\par\smallskip
If $a,b$ are two non-commuting elements of $p$-power order  then, by Lemma~\ref{lem-surjective}, there is an irreducible representation $W$ of $G$ such that $\bar a$ and $\bar b$ do not commute in $\bar G \subset \GL(W)$, and $\bar G$ is faithful of covariant dimension $2$. It follows from Proposition~\ref{prop-faithfulcovdim2} and Claim~1 that $\bar G/Z(\bar G) \simeq D_{2n}$. Since $p>2$, the images of $\bar a$ and $\bar b$ are in the cyclic subgroup $C_n\subset D_{2n}$ which leads to a contradiction, since the inverse image of $C_n$ in $\bar G$ is commutative.  
Hence pairs of elements of $p$-power order commute, and so the $p$-Sylow subgroup $G_p$ is normal and commutative.
It follows that $G':=\prod_{p>2} G_p$ is a normal commutative subgroup and that $G = G_2\cdot G'$ is a semidirect product. This proves Claim~2.

\smallskip
%
Now we can finish the proof. The case that $G=G_2$   is handled in  
Lemma~\ref{lem-2group}  below, so  
we can assume that   $G'$ is non-trivial. If $G_2$ commutes with $G'$, then $G_2$ is not commutative and faithful. Moreover, no $G_p$ can be of rank 2, else we have a subgroup which is a product $H:=G_2\times(\ZZ/p)^2$, and we have $\cdim H\geq 3$ by Corollary~\ref{cor-twoZ/p}. So $G'$ has rank  1. 
Then  $G'$ is cyclic, hence $G$ is faithful (Corollary \ref{cor-addZm}), which is a contradiction.  Hence we may assume that $G_2$ acts nontrivially on $G'$.

It is clear that $N_G=N_2\times N'$ where $N_2=N_G\cap G_2$ and $N':=N_G\cap G'$. Since $G_2$ acts nontrivially on $G'$, there is a $g\in G_2$ which induces an order 2 automorphism of some $G_p\neq\{e\}$. Then one can see that $g$ acts nontrivially on $N_{G_p}$. 
Since $G$ is not faithful, $N_G$ is not generated by a conjugacy class (Proposition \ref{prop-Gaschutz}) and the same holds for the subgroupÊ $H:=\langle g,N_2\rangle\ltimes N'$ (Corollary~\ref{cor-Gaschutz}). Thus $H$ is neither faithful nor commutative, so that it must equal $G$ by minimality. 

Suppose that $G_p=(\ZZ/p)^2$ for some $p$. If $g$ acts trivially on $G_p$, then it must act nontrivially on some $G_q$, and then we have the subgroup $(\langle g\rangle\ltimes G_q)\times (\ZZ/p)^2$ which by Corollary~\ref{cor-twoZ/p}  has covariant dimension at least $3$. 
If $g$ acts by sending each element of $(\ZZ/p)^2$ to its inverse, then $\langle g\rangle\ltimes G_p$ has covariant dimension 3 by Corollary~\ref{cor-Z2M}.
So we can assume that  $g$ acts on $G_p$ fixing one generator and sending the other to its inverse for every $G_p$ of rank 2. Thus $G'$ is generated by the conjugacy class of a single element. It follows that $N_2$ must have rank 2.  Moreover, $g$ must commute with $N_2$, else $N_2\times G'$ is generated by the conjugacy class of a single element. Suppose that $\langle g\rangle\cap N_2\simeq\ZZ/2$. If $g$ acts nontrivially on $\ZZ/p\subset G'$, then $\langle g,N_2\rangle\ltimes\ZZ/p$ contains a subgroup $(\langle g\rangle\ltimes\ZZ/p)\times\ZZ/2$ which has covariant dimension 3 by Corollary~\ref{cor-center}. If $\langle g\rangle\cap N_2=\{e\}$, then we have the subgroup $(\langle g\rangle \ltimes\ZZ/p)\times (\ZZ/2)^2$ which has covariant dimension three by Corollary~\ref{cor-twoZ/p}. This finishes the proof of the theorem, modulo the following lemma.
\end{proof}

\begin{lem}\label{lem-2group}
Let $G$ be a non-faithful $2$-group of covariant dimension 2. Then $G$ is commutative.
\end{lem}
\begin{proof}
Let $G$ be a counterexample of minimal order and let  $a,b\in G$ be two non-commuting elements. If $H:= \langle a,b \rangle$ is a strict subgroup of $G$, then $H$ is faithful and so $Z(H)$ has rank 1. Since $Z(G)$ has rank 2 this implies that $G$ contains a subgroup of the form $H\times\ZZ/2$ which has covariant dimension $3$  (Corollary~\ref{cor-center}). 

Thus we can assume that every pair $a,b\in G$ of non-commuting elements generates $G$. From the minimality assumption it follows that $a^2,b^2, (ab)^2 \in Z(G)$. Denote by $d:=(a,b)$ the commutator. Then $d = aba^{-1}b^{-1} = (ab)^2 a^{-2}b^{-2} \in Z(G)$.  
Since $aba^{-1}= db$ we have  
$b^2 = ab^2a^{-1} = d^2b^2$ which implies that $d^2 = 1$. Hence  $(G,G) = \{d,e\}$ and so $(a',b') = d$ for every pair $a',b'$ of non-commuting elements. Since $G/\langle d \rangle$ is commutative, it follows that $Z(G) = \langle a^2,b^2,d\rangle = \langle a^2,b^2,(ab)^2\rangle$. 

If  
$a^2 = z^2$ for some $z\in Z(G)$, then $a' := az^{-1} \notin Z(G)$ and ${a'}^2 = 1$. Thus $G$ contains a subgroup isomorphic to $Z(G)\times \ZZ/2$ which has covariant dimension 3. It follows that for every pair of non-commuting elements $a,b$  
the three elements $a^2$, $b^2$ and $(ab)^2$ are in $Z(G) \setminus Z(G)^2$. As a consequence, two of them generate $Z(G)$, i.e., there are generators $a,b$ of $G$ such that $a^2,b^2$ generate $Z(G)$. 
If $a^{2^p}=b^{2^q}\neq e$ where $p<q$ (and necessarily, $p$, $q\geq 2$), then the squares of $a':=ab^{-2^{q-p}}$ and $b$ freely generate $Z(G)$. We make a similar modification if $p>q$. If $p=q$, then we replace $a$ by $a':=ab\inv$. Then $(a')^{2^p}=a^{2^p}b^{-2^p}=e$. Thus we can assume that $a$ and $b$ freely generate $Z(G)$. This implies that if $a$ is of order $2^\ell$ and $b$ of order  $2^k$, then $\ell,k \geq 2$ and $|G| = 2^{\ell+k}$.

We have seen above that  $d^2 = 1$ where $d=(a,b)$. Therefore we are in one of the following three cases: (1) $d = a^{2^{k-1}}$, (2) $d=b^{2^{\ell -1}}$ or (3) $d=a^{2^{k-1}}b^{2^{\ell-1}}$.

\par\smallskip\noindent
{\bf Case 1:} $d = a^{2^{k-1}}$. Then $bab^{-1} = ad = a^{2^{k-1}+1}$ and so $G = \langle b \rangle \ltimes \langle a\rangle$  is a semidirect product.   
Then $\cdim G \geq 3$  by Proposition~\ref{prop-Z2Z2}.

\par\smallskip\noindent
{\bf Case 2:} $d = b^{2^{\ell - 1}}$. As in the previous case $G= \langle a \rangle \ltimes \langle b\rangle$  is a semidirect product and    
so  $\cdim G \geq 3$.

\par\smallskip\noindent
{\bf Case 3:} $d=a^{2^{k-1}}b^{2^{\ell-1}}$. We can assume that $\ell\geq k$. If $\ell>k$ then $a':=ab^{2^{\ell - k}}$ and $b$ generate $G$, $a'$ has the same order as $a$ and $d = {a'}^{2^{k-1}}$. Thus we are in Case~1 and so $\cdim G \geq 3$.

If $\ell=k>2$, we have $(ab)^2 = d a^2b^2$, and so $(ab)^{2^{\ell-1}} = d$. This implies that $G = \langle a \rangle \ltimes \langle ab \rangle$, hence $\cdim G \geq 3$.

Finally, for $\ell=k=2$ we get $(ab)^2 = e$ and so $G \supset \langle a^2,b^2,ab \rangle \simeq (\ZZ/2)^3$, hence again $\cdim G \geq 3$.
\end{proof}

Now we prepare the proof of Theorem \ref{thm-covdim2}. We need only consider faithful groups, and we can employ Proposition~\ref{prop-faithfulcovdim2}. So, we have an exact sequence
$$
1\to \ZZ/m\to G\to K\to 1
$$
where $K$ is $D_{2n}$, $n\geq 2$, $A_4$, $S_4$ or $A_5$. Thus we need to classify the cyclic central extensions of these groups. In terms of group cohomology, we need to calculate $\HH^2(K,\ZZ/m)$. First we determine the \emph{Schur multiplier} $\MM(K):=\HH^2(K,\CC^*)$.

\begin{lem}\label{lem-schurmult} We have
\begin{enumerate}
\item $\MM(A_5)\simeq\MM(S_4)\simeq \MM(A_4)\simeq\ZZ/2$.
\item $\MM(D_{2n})\simeq \ZZ/2$  if $n$ is even and $\MM(D_{2n})$  is trivial if $n$ is odd.
\end{enumerate}
\end{lem}
\begin{proof}
The first   part  is classical and goes back to Schur \cite{Schur}. The second part is surely also classical, but we don't know a reference, so we give a proof. 

Suppose that 
$$
1\to\CC^*\to H\to D_{2n}\to 1
$$
is exact. Then there are $\alpha$, $\beta_1\in H$ such that the image of $\alpha$ in $D_{2n}=\ZZ/2\ltimes\ZZ/n$ generates $\ZZ/2$ and the image of $\beta_1$ generates $\ZZ/n$. It is easy to arrange that $\alpha$ has order 2 and that $\beta_1$ has order $n$. Then $\alpha\beta_1\alpha\inv=\lambda\beta_1\inv$ where $\lambda\in\CC^*$ and $\lambda^n=1$. Replacing $\beta_1$ by $\beta_2:=\lambda^{-1/2}\beta_1$ we have that $\alpha\beta_2\alpha\inv=\beta_2\inv$ where the order of $\beta_2$ is now $n$ or $2n$. Suppose that the order is $2n$. If $n$ is odd, then $\beta:=\beta_2^2$ has order $n$ and maps to a generator of $\ZZ/n$. Thus if $\beta_2$ has order $n$ or $n$ is odd, our exact sequence is split. If $n$ is even, we see that there is a unique nontrivial extension, i.e., $\MM(D_{2n})\simeq\ZZ/2$. In fact, this extension is induced by  the nontrivial extension
$$
1\to \{\pm 1\}\to D_{4n}\to D_{2n}\to 1.
$$
\end{proof}

\begin{cor} \label{cor-H2} We have
\begin{enumerate}
\item $\HH^2(A_5,\ZZ/m)\simeq \ZZ/2$ if $m$ is even, else it is trivial.
\item $\HH^2(S_4,\ZZ/m)\simeq (\ZZ/2)^2$ if $m$ is even, else it is trivial.
\item $\HH^2(A_4,\ZZ/m)\simeq \ZZ/d$ where $d=\GCD(6,m)$.
\item $H:=\HH^2(D_{2n},\ZZ/m)\simeq (\ZZ/2)^3$ if $m$ and $n$ are even, $H\simeq \ZZ/2$ if $n$ is odd and $m$ is even and $H$ is trivial if $m$ is odd.
\end{enumerate}
\end{cor}
\begin{proof} We just give the proof of (c). The other proofs follow the same reasoning. From the short exact sequence
$$
1\to \ZZ/m\to\CC^*\xrightarrow{m}\CC^*\to 1 
$$
we obtain a long exact sequence of cohomology 
$$
\dots\Hom(A_4,\CC^*)\xrightarrow{m}\Hom(A_4,\CC^*)\to \HH^2(A_4,\ZZ/m)\to \MM(A_4)\xrightarrow{m}\MM(A_4)\dots
$$
where we use the fact that, since $\CC^*$ is a trivial $A_4$-module, we have  $\HH^1(A_4,\CC^*)\simeq\Hom(A_4,\CC^*)\simeq  \ZZ/3$. Now (c) follows from the exact sequence and Lemma \ref{lem-schurmult}.
\end{proof}

\begin{proof}[Proof of Theorem \ref{thm-covdim2}] We can assume that $G$ is faithful. Proposition~\ref{prop-faithfulcovdim2} gives us the possibilities for $G/Z(G)$. Suppose that $G/Z(G)\simeq  A_4$. If $\GCD(6,m)=1$, then we have a product extension of $A_4$ which has covariant dimension at least 
$\cdim A_4=3$. Suppose that we have a nonzero element of $\HH^2(A_4,\ZZ/m)$ which has order $3$. Then $3|m$ and we have the semidirect product $\ZZ/3m\ltimes(\ZZ/2)^2$ where the generator $\alpha$ of $\ZZ/3m$ permutes the nonzero elements of $(\ZZ/2)^2$ cyclically.   If $2|m$, then $G$ contains a copy of $(\ZZ/2)^3$, otherwise  $G\times\ZZ/2$ has the same covariant dimension as $G$ and contains a copy of $(\ZZ/2)^3$. Hence $\cdim G\geq 3$. Now suppose that we have an element of $\HH^2(A_4,\ZZ/m)$ of order $2$. Then we have the extension
$$
1\to\ZZ/m\to \widetilde {A_4}\to A_4\to 1
$$
where $\widetilde{A_4}\subset\GL_2$ is the binary tetrahedral group $BA_4$ multiplied by the $2m$th roots of 1 (as scalar matrices). So $\cdim G=2$. Finally, suppose that we have an element of order $6$. Here we need to be specific about the binary tetrahedral group. If we identify $\SU_2$ with the unit quaternions, then   $BA_4$ is generated by the subgroup $BK:=\{\pm i,\pm j,\pm k\}$ (the inverse image of the Klein 4-group) and the element $\alpha:=1/2(-1+i+j+k)$ which has order 3, so that $BA_4\simeq\ZZ/3\ltimes BK$. To get the extension of order $6$ we need to take the group generated by $BK$ and the product of $\alpha$ by a primitive $3m$th root of unity. But this is still a subgroup of $\GL_2$. This completes the proof for $A_4$.  

Now suppose that $G/Z(G)=S_4$. We can think of $S_4$ as the   semidirect product $\ZZ/2\ltimes A_4\simeq \ZZ/2\ltimes(\ZZ/3\ltimes(\ZZ/2)^2)$. If we have a trivial extension of $S_4$ by $\ZZ/m$, then we have a group of covariant dimension at least 3. If $m$ is even, then $\HH^2(S_4,\ZZ/m)$ has order 4, where the three nonzero elements correspond to the following three groups:
\begin{enumerate}
\item $\ZZ/2m\ltimes A_4$ where $\ZZ/2m$ has a generator which acts by an automorphism of order 2. This group contains $A_4$, hence has  covariant dimension at least 3. 
\item $(\ZZ/2m\ltimes BA_4)/(\ZZ/2)$ where the generator  $\alpha\in\ZZ/2m$ acts by order 2 on $BA_4$  and we identify $\alpha^m$ with $-1\in BA_4$. This is easily seen to be isomorphic to a subgroup of $\GL_2$. If $m=2$ it is the binary octahedral group $BS_4$.
\item $(\ZZ/m\times BS_4)/(\ZZ/2)$ where if $\alpha$ generates $\ZZ/m$ (acting as scalar matrices), then we identify $\alpha^{m/2}$ with $-1\in BS_4$. This  is again a subgroup of $\GL_2$.
\end{enumerate}

We now consider the case of $K:=D_{2n}$.  If the class in $\HH^2(K,\ZZ/m)$ is trivial, then we have $G=D_{2n}\times\ZZ/m$. If $n$ is odd then $D_{2n}$ has trivial center, so that $\cdim G=\cdim D_{2n}=2$ and $G$ is isomorphic to a subgroup of $\GL_2$. If $n$ is even, then we have center $\ZZ/2$ which means that we have covariant dimension $3$ if $m$ is even and covariant dimension 2 (and $G$ is a subgroup of $\GL_2$) if $m$ is odd. From now on we can suppose that $m$ is even. 

Suppose that $n$ is odd and that we have a nontrivial extension of $D_{2n}$. Then the only candidate is
$(\ZZ/m\times D_{4n})/\ZZ/2$ where we identify the $(m/2)$nd power of the generator of $\ZZ/m$ with the central element in $D_{4n}$. This is clearly a subgroup of $\GL_2$. From now on we can assume that $n$ is even.

Choose $\alpha\in G$ whose   image generates $\ZZ/2\subset D_{2n}$ and $\beta\in G$ whose image generates $\ZZ/n\subset D_{2n}$ where $D_{2n}=\ZZ/2\ltimes\ZZ/n$. Then $\alpha^2\in\ZZ/m$ and $\beta^n\in\ZZ/m$. We have that $\alpha\beta\alpha\inv=z\beta\inv$ where $z\in\ZZ/m$. Replacing $\beta$ by a product $\beta z'$ for $z'\in\ZZ/m$ we can reduce to the case that $z=e$ or that $z$ is a (fixed) generator of $\ZZ/m$. Similarly, we can assume that $\alpha^2=w$ where  $w=e$ or $w$ is the same fixed generator of $\ZZ/m$ as above. Now $\alpha\beta^n\alpha\inv=z^n(\beta\inv)^n=\beta^n$, so that $\beta^{2n}=z^n$ and $\beta^n=\pm z^{n/2}$ (where we think of $-1$ as $n/2\in\ZZ/n$). If $\beta^n= z^{n/2}$, then it follows that $\alpha$ fixes $\beta^{n/2}$, hence that $\beta^{n/2}\in\ZZ/m$ which would imply that the image of $\beta$ in $\ZZ/n$ would have order $n/2$, a contradiction. Hence $\beta^n=-z^{n/2}$.

\par\smallskip\noindent 
{\bf Case 1:} $\alpha^2=e$. If $z=e$, then $\beta$ has order $n$ so that our exact sequence is split, a case that we have already handled. So we can assume that $\beta^n=-z^{n/2}$ where $z$ generates $\ZZ/m$. Now we see that $G\subset\GL_2$. Represent $\alpha$ by the matrix $\left(\smallmatrix 0 & 1 \\ 1 & 0 \endsmallmatrix\right)$ and represent $\beta$ as the matrix $\left(\smallmatrix \xi\eta & 0 \\ 0 & \xi\inv\eta \endsmallmatrix\right)$ where $\xi$ is a primitive $2n$th root of $1$ and $\eta$ is a primitive $2m$th root of 1. Then $\beta^n$ is central and we have that $\alpha\beta\alpha\inv=z\beta\inv$ where $z=\left(\smallmatrix \eta^2 & 0 \\ 0 & \eta^2 \endsmallmatrix\right)$ generates a central copy of $\ZZ/m$  and $\beta^n=-z^{n/2}$. 

\par\smallskip\noindent
{\bf Case 2:} $\alpha^2=w$ generates $\ZZ/m$. If $z=e$, then our group is %
isomorphic to the group
generated by the matrices $\left(\smallmatrix 0&  \eta \\   \eta & 0  \endsmallmatrix\right)$ where $\eta$ is a primitive $2m$th root of 1 and 
$\left(\smallmatrix \xi  & 0 \\ 0 & \xi\inv  \endsmallmatrix\right)$ where $\xi$ is a primitive $2n$th root of 1. If $z=w$ generates $\ZZ/m$, then we use $\alpha$ as above and $\beta= \left(\smallmatrix \xi\eta & 0 \\ 0 & \xi\inv\eta \endsmallmatrix\right)$ where $\xi$ is a primitive $2n$th root of 1. Then $\alpha$ and $\beta$ generate %
a subgroup of $GL_2$ isomorphic to $G$.
\end{proof}
\vskip1cm

\end{document}